\documentclass[12pt]{article}
\usepackage{amsmath,amsfonts,amssymb,latexsym,amsbsy, bbm, theorem,enumerate,color}
\usepackage{graphicx, graphics}
\usepackage{float,color,fancybox,shapepar,setspace,hyperref}
\usepackage[affil-it]{authblk}
\usepackage{caption}
\usepackage{tikz}
\usetikzlibrary{arrows}
\usepackage[normalem]{ulem}
\usetikzlibrary{arrows,automata}
\usepackage[latin1]{inputenc}
\usepackage{verbatim}

\usepackage{subfigure}
\usepackage{psfrag}
\usepackage{epstopdf}

\newtheorem{thm}{Theorem}[section]

\newtheorem{prop}[thm]{Proposition}
\theorembodyfont{\rmfamily}
\newtheorem{defn}[thm]{Definition}

\begin{document}

\title{The average eccentricity of a graph with prescribed girth}

\author{Fadekemi Janet Osaye\thanks{Some of the first result presented in this paper forms part of the  author's PhD thesis at the University of Johannesburg, South Africa. Financial support from 2018 Kovaleskaia Research Grants for Female Mathematicians is gratefully acknowledged.}}

  \affil{
  { \small {Department  of Mathematics and Statistics, Auburn University, AL  36830, AL}}\\
  }


\maketitle
\begin{abstract}
	Let $G$ be a connected graph of order $n$. The eccentricity $e(v)$ of a vertex $v$ is the distance from $v$ to a vertex farthest from $v$. The average eccentricity of $G$ is the mean of all eccentricities in $G$. We give upper bounds on the average eccentricity of $G$ in terms of order $n$, minimum degree $\delta$, and girth $g$. In addition, we construct graphs to show that, if for given $g$ and $\delta$, there exists a Moore graph of minimum degree $\delta$ and girth $g$, then the bounds are asymptotically sharp. Moreover, we show that the bounds can be improved for a graph of large degree $\Delta$.

\end{abstract}

\baselineskip 18pt
\section{Introduction}

Let $G$ be a connected graph with vertex set $V(G)$. The {\em eccentricity} $e(v)$ of a vertex is the maximum distance between $v$ and a vertex in $G$. 
The {\em average eccentricity} ${\rm avec}(G)$ of $G$ is the mean of all eccentricities of the vertices of $G$, i.e., ${\rm avec}(G) = \frac{1}{|V(G)|}\sum_{u\in V(G)} e(u)$. The average eccentricity, introduced by Buckley and Harary as \textit{eccentric mean} \cite{BucHar}, was originally conceived as a performance indicator for transportation networks. It later attracted much interests within graph theory. Several results have also been inspired by conjectures of the computer programme AutoGraphix by Aouchiche et. al.\cite{Aouchiche}. The first systematic study on average eccentricity was initiated by Dankelmann, Goddard and Swart \cite{DanGodSwa2004}, who, amongst other results showed that for all connected graphs of order $n$,
\begin{equation}\label{eqn1a}
{\rm avec}(G) \leq \frac{1}{n}\left\lfloor \frac{3n^2}{4} -\frac{n}{2}\right\rfloor,
\end{equation}
with equality if and only if $G$ is a path of order $n$. \\

Several upper bounds on the average eccentricity in terms of order and minimum degree are known. Dankelmann, Goddard and Swart \cite{DanGodSwa2004} proved that for all connected graphs $G$ of order $n$ and minimum degree $\delta \geq 2$,

\begin{equation}\label{eqn1}
{\rm avec}(G)  \leq \frac{9n}{4(\delta +1)} +\frac{15}{4}, 
\end{equation}
and this bound is sharp apart from a small additive constant. By using similar argument initiated by Dankelmann and Entringer \cite{DanEnt2000}, stronger bounds for triangle-free graphs and $C_4$-free graphs \cite{DanMukOsaRod} were proven, respectively, to be \begin{equation}\label{eqn2}
{\rm avec}(G)  \leq 3\left\lceil \frac{n}{2\delta}\right\rceil +5, 
\end{equation}
and, 
\begin{equation}\label{eqn3}
{\rm avec}(G) \leq \frac{15}{4} \left\lceil \frac{n}{\delta^2 -2\lfloor \frac{\delta}{2} \rfloor +1} \right\rceil +\frac{11}{2}. 
\end{equation}
Also, it was shown in \cite{AlexDan2020} that if $G$ is a graph of girth at least 6, then 
\begin{equation}\label{eqn4}
{\rm avec}(G)  \leq \frac{9}{2} \left\lceil \frac{n}{2\delta^2 -2\delta +2} \right\rceil + 8,
\end{equation}
which when relaxed for $(C_4,C_5)$-free graphs yields

\begin{equation}\label{eqn5}
{\rm avec}(G)  \leq \frac{9}{2} \left\lceil \frac{n}{2\delta^2 -5\delta +5} \right\rceil + 8.
\end{equation}

%

However, stronger bounds obtained by Dankelmann and Osaye \cite{DanOsa2019} showed that for fixed $\delta$ and maximum degree $\Delta$, inequations \eqref{eqn1}, \eqref{eqn2} and \eqref{eqn3} can be improved, respectively, to 

\begin{equation}\label{eqn6}
{\rm avec}(G)  \leq \frac{9(n-\Delta-1)}{4(\delta +1)}\left(1+ \frac{\Delta -\delta}{3n}\right) +7
\end{equation}
\begin{equation}\label{eqn7}
{\rm avec}(G)  \leq \frac{3(n-\Delta)}{2\delta}\left(1+\frac{\Delta -\delta}{3n}\right) +\frac{19}{2}, 
\end{equation}
and 

\begin{equation}\label{eqn8}
{\rm avec}(G) \leq
\frac{15}{4} \frac{n-\varepsilon_{\Delta} + \varepsilon_{\delta}}{\varepsilon_{\delta}} 
\Big[ 1 + \frac{\varepsilon_{\Delta} - \varepsilon_{\delta}}{3n} \Big]  + \frac{37}{4},
\end{equation}
where 
$\varepsilon_{\Delta}:= \Delta \delta -2\Big\lfloor\frac{\Delta}{2}\Big\rfloor +1$
and 
$\varepsilon_{\delta} := \delta^2 -2\lfloor \frac{\delta}{2} \rfloor +1$. 
\\
Several other bounds or relations on the average eccentricity with given properties or known parameters can be found in the literature. For example, for trees of given  maximum degree \cite{DuIli2013}, for graphs of given order and size (\cite{AliDanMorMukSwaVet2018}, \cite{TanZho2012}), order and independence or domination number \cite{DanMuk2014}, and for graphs of given order and $k$-packing or $k$-domination number \cite{DanOsa-manu}. 
Remarkably, very little is known about generalised forms for distance parameters. For example, Ali \cite{PAli} showed that the previous bounds on the Steiner diameter of a graph can be generalised when the girth of the graph is prescribed. This aim of this paper is to investigate the previous results on the average eccentricity of a graph and provide a generalised bound given its girth. Using methods similar to the one used in \cite{DanMukOsaRod} and \cite{DanOsa-manu}, we show that, for connected graphs $G$ of given $\delta \geq 3$ and girth $g$ odd with $K= 1 + \frac{\delta}{\delta-2}\left[(\delta -1)^{\frac{g-1}{2}} -1 \right] $,  
\[{\rm avec}(G)  \leq \frac{3g}{4}\left\lceil \frac{n}{K}\right\rceil+\frac{3g}{2}-2.\] For $g$ even with $L= \frac{2}{\delta-2}\left[(\delta -1)^{\frac{g}{2}} -1 \right] $, 

\[{\rm avec}(G)  \leq \frac{3g}{4}\left\lceil \frac{n}{L}\right\rceil+\frac{3g}{2}-2.\] 
In addition, we construct graphs to show that, if for given $\delta$ and $g$, there exists a Moore graph of minimum degree $\delta$ and girth $g$, then the bounds are asymptotically sharp. Moreover, we obtain stronger bounds for graphs of fixed $\delta$ and maximum degree $\Delta$. These results are, in fact, generalisations of inequations \eqref{eqn1} - \eqref{eqn8}. \\
\\
The notation used is as follows. For a connected graph $G$ of order $n$ with vertex set $V(G)$ and edge set $E(G)$, the distance between two vertices $u$ and  $v$, $d_G(u,v)$, in $G$ is the length of a shortest $u-v$ path in $G$. The eccentricity of a vertex $v$, $e_G(v)$, is the distance between $v$ and a vertex farthest from $v$ in $G$. The diameter ${\rm diam}(G)$ and ${\rm rad}(G)$ are, respectively, the largest and minimum eccentricity of vertices of $G$. If no confusion arise, we will drop the subscript  $G$ throughout the paper. The total eccentricity, $EX(G)$ is the sum of all eccentricities in $G$. The open neighbourhood of a vertex $v$, $N(v)$, is the set of vertices adjacent to $v$. The closed neighbourhood of $v$, $N[v]$, is defined as $N(v) \cup \{v\}$. For $A\in V(G)$ and $k\in \mathbb{N}$, the $k$-th neighbourhood of $A$, $N^k(A)$, is the set of all vertices $x$ of $G$ with $d(x,a)\leq k$ for some $a\in A$. The $k$-th power of $G$, denoted as $G^k$, is the graph with the same vertex set as $G$, in which any two vertices $u, v \in V(G)$ are adjacent if $d(u,v)\leq k$. For a positive integer $k$, the $k$-packing of $G$ is a subset $A \subset V(G)$ with $d(a,b) > k$, for all $a,b \in A$. If $B \subset V(G)$, then $G^k[B]$ is the subgraph induced by $B$.\\
\\
If $H$ is a subgraph of $G$, we write $H\leq G$. A set $M\subset E(G)$ in $G$ is a matching in $G$ if no two edges of $M$ are incident. A matching of maximum size is called a maximum matching. The line graph $L(G)$ of $G$ is the graph whose vertices are the edges of $G$ such that two vertices of $L(G)$ are adjacent if they share a vertex as edges of $G$. Let $M\subseteq E(G)$, then $V(M)$ is the set of vertices incident with at least an edge of $M$.
\section{Results}
We will need the following definition and proposition due to Dankelmann, Goddard and Swart. \cite{DanGodSwa2004}. 

\begin{defn}[\cite{DanGodSwa2004}]\label{def1}
	Let $G$ be a connected graph and $c:V(G) \rightarrow  \mathbb{R}$ a 
	nonnegative weight function
	on the vertices of $G$. Then the eccentricity of $G$ with respect to $c$ is defined by 
	\[ EX_c(G)= \sum_{x\in V(G)} c(x) e_G(x). \] 
	Let $N=\sum_{x\in V(G)} c(x)$ be the total weight of the vertices in $G$. If $N>0$, then the average eccentricity of $G$ with respect to $c$ is 
	\[ {\rm avec}_c(G) = \frac{EX_c(G)}{N}.   \]
\end{defn}

\begin{prop}[\cite{DanGodSwa2004}]\label{prop2}
	Let $G$ be a connected graph, $c$ a weight function on the vertices 
	of $G$, and $N = \sum_{v\in V(G)} c(v)$ the total weight of the vertices
	of $G$. If $c(v) \geq 1$ for all $v\in V(G)$, then 
	\[
	{\rm avec}_c(G) \leq {\rm avec}(P_{\lceil N \rceil}).
	\]
\end{prop}
\subsection{Bounds in terms of order, minimum degree and girth}
We now present one of our main results. 

\begin{thm}\label{theoremgirth}
	Let $G$ be a connected graph with $n$  vertices, minimum degree $\delta \geq 3$ and girth $g$.\\ 
	(a)  If $g$ is odd, then $${\rm avec} (G) \leq  \frac{3g}{4}\left\lceil \frac{n}{K}\right\rceil+\frac{3g}{2}-2,$$ where  $K= 1+ \frac{\delta}{\delta -2} [(\delta-1)^{\frac{g-1}{2}} -1]$.
	\\
	\\
	(b)  If $g$ is even, then  $${\rm avec} (G) \leq  \frac{3g}{4}\left\lceil \frac{n}{L}\right\rceil+\frac{3g}{2}-2,$$ where 
	$ L= \frac{2}{\delta-2}\Big[(\delta-1)^{\frac{g}{2}}-1\Big].$
	The bounds are sharp apart from a small additive constant.
\end{thm}
%

\noindent
\textbf{Proof.}
(a)\,\, We first find a maximal $(g-1)$-packing $A$ of $G$ using the following procedure. Choose an arbitrary vertex $a_1$ of $G$ and let $A= \{a_1\}$. If there exists a vertex $a_2$ in $G$ with $d_G(a_1, a_2)=g$, add $a_2$ to $A$. Add vertices with  $d_G(u, A)=g$ to $A$ until each of the vertices not in $A$ is within distance $(g-1)$ of $A$.\\
Let $T_1 \leq G$ be the subforest of $G$ with vertex set ${N}_G^{\frac{g-1}{2}}[A]$ in which for every $a\in A$, the set $N_G^{(g-1)/2}(a)$ forms
a component of $T_1$ so that the distances to $a$ are preserved. With the way $A$ is constructed, there exist $|A|-1$ edges in $G$, each joining two vertices of distinct components of $T_1$ whose addition to $T_1$ yields a subtree $T_2 \leq G$.\\
Now each vertex $u \in V(G) - V(T_2)$ is within distance $(g-1)/2$ to some vertex $v\in V(T_2)$. Let $T$ be the spanning tree of $G$ containing $T_2$ such that the distance of every vertex to $A$ is preserved.\\
Since deleting edges does not decrease the average eccentricity, ${\rm avec}(G) \leq {\rm avec}(T)$. Thus it suffices to show that 
\begin{equation*}\label{eqn2o}
{\rm avec} (T) \leq  \frac{3g}{4}\left\lceil \frac{n}{K}\right\rceil+\frac{3g}{2}-2.
\end{equation*}
For every  vertex $u\in V(T)$, let $u_A$ be the unique vertex in $A$ closest to $u$ in $T$. If we consider ${\rm avec}(T)$ as the weighted average eccentricity of $T$ under the weight function which is identically $1$ on all vertices, then the weight function $c$, which we define now, can be thought of as obtained by moving the weight of every vertex to the closest vertex in $A$. That is, we define a weight function $c:V(T) \rightarrow \mathbb{R}^+$ by $$c(u) = \Big|\{ x \in V(T)|  x_A=u\}\Big|, \, \text{for each}\,\, u\in A.$$ 
Since $G$ has girth $g$ and $g$ is odd, and all vertices in $G$ have at least degree $\delta$, we have 
\begin{equation}
\begin{aligned}
c(u) &\geq |{N}_G^{(g-1)/2}[u]|\\
&\geq 1+ \delta +\delta (\delta-1) + \delta (\delta-1)^2+ \ldots + \delta(\delta-1)^{(g-3)/2}.\\
&= 1+ \frac{\delta}{\delta-2}\big[ (\delta-1)^{(g-1)/2}-1\big] = K.\\
\end{aligned}
\end{equation}
We know that for every vertex $u$ of $T$, $d(u,u_A)\leq g-1$, thus $|e_T(u)-e_T(x_A)|\leq g-1$. Hence,
\begin{eqnarray}\label{eqn2p}
{\rm avec}(T) & = & \frac{1}{n} \sum_{x \in V(T)} e_T(u) \nonumber \\
& \leq & \frac{1}{n} \sum_{u\in V(T)} [e_T(u_A) + g-1] \nonumber \\
& = & \big[ \frac{1}{n} \sum_{u \in A} c(u) e_T(u) \big] + g-1 \nonumber \\
& = & {\rm avec}_{c}(T) + g-1.    \label{eqn-5b} 
\end{eqnarray}
Now the weight is concentrated on the vertices of $A$. By our construction of $A$, $T^g[A]$ is connected. Therefore, for all $a_i, a_j \in A$, $d(a_i,a_j) \leq g\,d_{T^g[A]}(a_i,a_j)$. Since each vertex of $T$ is within distance $g$ of some vertex in $A$, we have for each $u\in A$,  
$e_T(u) = g e_{T^g[A]}(u) + (g-1)$, which implies that
\begin{equation}\label{eqn2q}
{\rm avec}_c(T) \leq g\, {\rm avec}_c(T^g[A]) +(g-1).
\end{equation}
We define a new weight function by $c'(u) = \frac{c(u)}{K}$, and let $N'=  \sum_{u\in A} c'(u)=   \frac{n}{K} $. Hence, 
\begin{equation*}
\begin{aligned}
{\rm avec}_{c'}(T^g[A])& = \frac{EX_{c'}(T^g[A])}{\sum_{u\in A}c'(u)}\\
&= \frac{\frac{1}{\delta +1}\sum_{u\in A} c(u) e_{(T^g[A])}(u)}{\frac{1}{\delta+1} \sum_{u\in A}c(u)}\\
&= {\rm avec}_c (T^g[A]).
\end{aligned}
\end{equation*}
By Proposition \ref{prop2} and \eqref{eqn1a},  we have
\begin{equation}\label{eqn2r}
{\rm avec}_c(T^g[A]) \leq {\rm avec} (P_{\lceil N'\rceil}) \leq \frac{3\lceil N'\rceil}{4}-\frac{1}{2}.
\end{equation} 
Then by combining \eqref{eqn2p}, \eqref{eqn2q} and \eqref{eqn2r}, we have
\begin{equation*}
\begin{aligned}
{\rm avec}(T) &\leq  {\rm avec}_c(T) + g-1\\
& \leq g\, avec_c(T^g[A]) +2(g-1)\\
& \leq g\left(\frac{3\lceil N'\rceil}{4}-\frac{1}{2}\right)+2(g-1) \\
&=\frac{3g\lceil N'\rceil}{4}+ \frac{3g}{2}-2\\
&= \frac{3g}{4}\left\lceil \frac{n}{K}\right\rceil+\frac{3g}{2}-2.
\end{aligned}
\end{equation*}
This proves the (a) part of the theorem.\\
\\
(b)\,\, For $g$ even, we first find a maximal matching $M$ of $G$ using the following procedure. Choose an arbitrary edge $e_1 \in E(G)$ and let $M= \{e_1\}$. If there exist an edge $e_2$ in $G$ with $d_G(e_2,e_1)=g-1$, add $e_2$ to $M$. Add edges with $d_G(e_i,V(M))=g-1$ to $M$ until each of the edges not in $M$ is within distance at most $g-2$ of $M$. \\
\\
Let $T_1 \leq G$ be the forest subgraph with vertex set ${N}_G^{\frac{g-2}{2}} [M]$ in which for every $v\in V(M)$, the set $N_G^{(g-2)/2}(v)$ forms
a component of $T_1$ so that the distances to $v$ are preserved. By our construction of $M$, there exist $|M|-1$ edges in $G$, each joining two distinct components of $T_1$, whose addition to $T_1$ yields a tree $T_2 \leq G$. \\
\\
Now each vertex in $V(G)- V(T_2)$ is within distance $\frac{g-2}{2}$ of some vertex $u' \in V(T_2)$. Let $T$ be a spanning tree of $G$ containing $T_2$, such that the distance of every vertex of $T$ to $M$ is preserved.
Since ${\rm avec}(G) \leq {\rm avec}(T)$, it suffices to prove the result for $T$. 
For every $u\in V(T)$, let $u_M$ be a vertex in $V(M)$ closest to $u$ in $T$. As in the proof of part (a), we think of ${\rm avec}(T)$ as the average eccentricity with respect to the weight function which is constant $1$, and we move the weight of $u$ to $u_M$. That is, we define a weight function $c: V(T) \rightarrow \mathbb{R}^+$  by: $$ c(u)= \big|\{x\in V(T)| x_M = u \}\big|, \,\, \text{for all}\,\, u\in V(T).$$
Thus, for each $u\in M$,
\begin{equation*}
\begin{aligned}
c(u) &\geq |{N}_G^{(g-2)/2}[u]|\\
&\geq 1+ (\delta-1) + (\delta-1)^2+ \ldots + (\delta-1)^{(g-2)/2}\\
&= \frac{(\delta-1)^{g/2}-1}{\delta -2}.
\end{aligned}
\end{equation*}
Now the weight $c$ is concentrated exclusively on the vertices of $V(M)$. Consider the line graph $L=L(T)$ defined by the weight function $\bar{c}$ with $V(T)=E(T)$ such that:
\[
\bar{c}(uv) =
\begin{cases}
c(u) +c(v) & \text{if } uv \in M, \\
0
& \text{if } uv \not \in M.
\end{cases}
\] 
Then, $$\bar{c}(uv) \geq \frac{2(\delta-1)^{g/2}-2}{\delta -2}=L, \,\, \text{for each}\, uv \in M.$$
Since $d_T(u,v)=g-1$ for any two vertices $u,v\in V(M)$, each weight was moved over a maximum distance of $g-2$ to the nearest vertex in $V(M)$. Hence, 
\begin{equation}\label{2.20}
{\rm avec}(T) \leq {\rm avec}_c(T) + (g-2).
\end{equation}
If $e_1, e_2 \in E(T)$ are edges incident with vertices $v_1,v_2 \in V(T)$, respectively, and $v_1$ is an eccentric vertex of $v_2$ in $T$, then $e_T(v_2) = d_T(v_1,v_2) \leq d_L(e_1,e_2) +1 \leq e_L(e_2) + 1$. It follows that

\begin{equation}\label{2.21}
{\rm avec}_c(T) \leq {\rm avec}_{\bar{c}}(L) +1.
\end{equation}
If the distance $d_T(e_1,e_2)$ between any two matching edges $e_1,e_2 \in M$ equals $(g-1)$, then $d_L(e_1,e_2) \leq g$. Hence  $L^g[M]$ is a connected graph. Thus, 
\begin{equation}\label{2.22}
{\rm avec}_{\bar{c}}(L) \leq g\, {\rm avec}_{\bar{c}}(L^g[M]) + g-1.
\end{equation}
Now let $\bar{c}'(uv)= \frac{\bar{c}(uv)}{L}$ be a new weight function, and let $N'= \sum_{uv\in M}\bar{c}'(uv)$, for each $uv\in M$. Then $N'= \frac{n}{L}$. As in the case of (a) above, $${\rm avec}_{\bar{c}}(L^g[M])= {\rm avec}_{\bar{c'}}(L^g[M]).$$ Since $L^g[M]$ is connected,  we have by Proposition \ref{prop2} and equation \eqref{eqn1a}, 
\begin{equation}\label{2.23}
{\rm avec}_{\bar{c}}(L^g[M])\leq {\rm avec}(P_{N'})\leq  \frac{3\lceil N'\rceil}{4}- \frac{1}{2}.
\end{equation}
Thus by substituting \eqref{2.20}, \eqref{2.21} and \eqref{2.22} in \eqref{2.23}, we have 
\begin{equation*}
\begin{aligned}
{\rm avec}(T) &\leq {\rm avec}_c(T) + g-2\\
& \leq g\, {\rm avec}_{\bar{c}}(L) +g-1\\
& \leq g\, {\rm avec}_{\bar{c}}(L^g[M]) +2(g-1)\\
& \leq g\left(\frac{3\lceil N'\rceil }{4}-\frac{1}{2}\right)+2(g-1) =\frac{3g\lceil N'\rceil }{4}+ \frac{3g}{2}-2.\\
\end{aligned}
\end{equation*}
The (b) part of the theorem now follows since $N'=\frac{n}{L}$. 
\\
\\
\\
\noindent
We verify that these bounds are true for previous results with given order and minimum degree or with maximum degree as shown in Equations \eqref{eqn1} to \eqref{eqn8}.
\begin{itemize}
	\item[(a)] If $g=3$, then $K=\delta +1$, so by a simple calculation in conjunction with $\left\lceil \frac{n}{\delta +1} \right\rceil < \frac{n}{\delta +1} +1$, ${\rm avec}(G) \leq \frac{9}{4} \left\lceil \frac{n}{\delta +1} \right\rceil +\frac{5}{2} < \frac{9n}{4(\delta +1)}+ \frac{15}{4}$. Note that this is equivalent to the bound \eqref{eqn1} obtained in \cite{DanGodSwa2004}. 
	\item [(b)] If $g=4$, then $G$ is triangle-free and $L =2\delta$. Thus, by a simple calculation, ${\rm avec}(G)\leq 3\left\lceil \frac{n}{2\delta}\right\rceil +4$, which differs from \eqref{eqn2} proved in \cite{DanMukOsaRod} by an additive constant 1.   
	\item[(c)]  If $g=5$, then $G$ is $(C_3,C_4)$-free and $K= \delta^2 +1$. By substituting $g$ and $L$ in part $(b)$ of the theorem, we have ${\rm avec}(G) \leq \frac{15}{4} \left\lceil \frac{n}{\delta^2 +1} \right\rceil + \frac{11}{2}$, which is a slightly weaker bound to \eqref{eqn3} proved in \cite{DanMukOsaRod} for strictly $C_4$-free graphs.
	\item[(d)] If $g=6$, then $G$ is $(C_3, C_4, C_5)$-free and $L=\frac{2}{\delta -2} \left( (\delta -1)^3-1\right)$. A simple calculation shows that ${\rm avec}(G) \leq \frac{9}{2} \left\lceil \frac{n(\delta -2)}{2[(\delta -1)^3-1]} \right\rceil + 7 \leq \frac{9}{2} \left\lceil \frac{n}{2(\delta^2-\delta +1)} \right\rceil + 7 $. This differs from \eqref{eqn4} proved in \cite{AlexDan2020} for graphs of girth at least 6, by an additive constant 1.
\end{itemize}
The bounds in Theorem \ref{theoremgirth} are asymptotically sharp, apart from small additive constants, if $\delta$ and $g$ are such that there exists a Moore graph, i.e., a graph with minimum degree $\delta$, girth $g$, diameter  $d= \frac{g-1}{2}$, and  order $ K = 1 + \frac{\delta}{\delta-2}\Big[(\delta-1)^{\frac{g-1}{2}}-1\Big]$ (for $g$ odd) or $ L= \frac{2}{\delta-2}\Big[(\delta-1)^{\frac{g}{2}}-1\Big]$ (for  $g$ even).\\
\\
For a given integer $k>0$, let $G_1, G_2, \ldots, G_k$ be disjoint copies of the $(\delta, g)$- Moore graph. Let $G_{n,\delta,k}$ be the graph obtained from the union of $G_1, G_2, \ldots, G_k$ by deleting the edges $a_ib_i$ for $i= 2,3, \ldots, k-1$ and adding the edges $a_{i+1}b_i$ for $i= 1,2, \ldots, k-1$. 
\begin{itemize}
	\item[(a)]  If $g$ is odd, then each graph $G_i-a_ib_i$ has diameter at least $2d$ since $d_{G_i-a_ib_i}(a_i,b_i) \geq 2d$. Since ${\rm diam}(G_i)=d$, for $i=1,k$, and since there are $k-1$ link edges of the form $a_{i+1}b_i$, we have in conjunction with $g=2d+1$, ${\rm diam}(G_{n,\delta,k})=g(k-1)$. 
	It is easy to see that ${\rm rad}(G_{n,\delta, k}) = \left\lceil \frac{g(k-1)}{2} \right\rceil$.
	For each vertex $x\in V(G_i)$, $1 \leq i \leq \lceil\frac{k-1}{2}\rceil$, 
	$$e_{G_i}(x) \geq d(b_i,b_k)=g(k-i) -d.$$ 
	Moreover, $\sum_{i=1}^{\lceil\frac{k-1}{2}\rceil}(G_i) = \sum_{i=1}^{\lceil\frac{k-1}{2}\rceil}(G_{k-i+1})$ for each $i$, $1 \leq i \leq \lceil\frac{k-1}{2}\rceil$.  Since each $G_i$ has $K$ vertices, $\sum_{i=1}^{\lceil\frac{k-1}{2}\rceil}(G_i) \geq K[g(k-i) -d]$.
	\begin{eqnarray*}
		EX(G_{n,\delta,k}) &=& 2\Big[\sum_{i=1}^{\lceil\frac{k-1}{2}\rceil} \sum_{x \in G_i}e_{G_{n,\delta,k}}(x)  \Big]\\
		&\geq & 2K [\sum_{i=1}^{\frac{k}{2}} g(k-i) -d]= 2K\Big[\sum_{i=1}^{\frac{k}{2}}g(k-i) -d(\frac{k}{2})\Big]\\
		&=& 2Kg\Big[\frac{3k^2}{8}-\frac{k}{4}\Big]-Kdk
	\end{eqnarray*}
Therefore,\begin{equation}\label{eqn2t}
{\rm avec}(G_{n,\delta,k})\geq \frac{gK}{n}\Big[\frac{3k^2}{4}-\frac{k}{2}\Big]-\frac{Kdk}{n}.
\end{equation}
Now \eqref{eqn2t} in conjunction with $K=\frac{n}{k}$ and $d=\frac{g-1}{2}$ yields
\[{\rm avec}(G_{n,\delta,k})\geq \frac{3gn}{4K}-g+\frac{1}{2}.\]
In this case, the bound between ${\rm avec}(G_{n,\delta,k})$ and the bound in Theorem \ref{theoremgirth}(a) differs only by an additive constant at most $\frac{5}{2}(g-1)$.  
\item[(b)] If $g$ is even, a simple calculation shows that ${\rm diam}(G_{n,\delta,k})=g(k-1)+1$ and ${\rm rad}(G_{n,\delta,k})=\frac{g(k-1)}{2}+1$; then for each vertex $x\in V(G_i)$, $1\leq i \leq \lceil\frac{k}{2}\rceil$, 
$e_{G_i}(x) \geq g(k-i) -d+1$.
For $g$ and $k$ even, $\sum_{i=1}^{\lceil\frac{k-1}{2}\rceil}(G_i)) \geq L[g(k-i)-d+1]$. Hence,
\begin{eqnarray*}
	EX(G_{n,\delta,k}) &=&2\Big[\sum_{i=1}^{\lceil\frac{k}{2}\rceil} \sum_{x \in G_i}e_{G_{n,\delta,k}}(x)  \Big]\\
	&\geq &  2L\Big[\sum_{i=1}^{\frac{k}{2}}g(k-i) -(d-1)(\frac{k}{2})\Big]\\
	&=& 2Lg\Big[\frac{3k^2}{8}-\frac{k}{4}\Big]-Lk(d-1).
\end{eqnarray*}
In conjunction with $L=\frac{n}{k}$ and $d=\frac{g-1}{2}$, we obtain as in the case of $g$ odd,
\[{\rm avec}(G_{n,\delta,k})\geq \frac{3gn}{4L}-g+\frac{3}{2}.\] In this case, this bound differs from the bound in Theorem \ref{theoremgirth} (b) by at most $\frac{5}{2} (g-1) -1$.
\end{itemize}
\hfill $\Box$
\\
\\

\subsection{Bounds in terms of order, girth, minimum degree and maximum degree}  
In Theorem \ref{theoremgirth}, the examples that prove that the bounds are asymptotically sharp proposed that all vertices have degrees close to the minimum degree. This suggests that we can improve these bounds if $G$ contains a vertex of large degree.
We now presents bounds on the average eccentricity of a graph that contains a vertex of maximum degree $\Delta$. The proofs follow closely the methods used in Theorem  \ref{theoremgirth} and \cite{DanOsa2019} with some slight modifications.

\begin{thm}\label{theoremgirth2}
	Let $G$ be a connected graph with $n$ vertices, minimum degree $\delta$, maximum degree $\Delta$ and girth $g$.\\ 
	(a)  If $g$ is odd, then $${\rm avec}(G)\leq   \frac{3g}{4} \Big(\frac{n-K_2}{K_1} \Big) \Big(1 + \frac{K_2-K_1}{3n} \Big) 
	+ 3g -2,$$ where $ K_1 = 1 + \frac{\delta}{\delta-2}\Big[(\delta-1)^{\frac{g-1}{2}}-1\Big]$ and $ K_2 = 1 + \frac{\Delta}{\delta-2}\Big[(\delta-1)^{\frac{g-1}{2}}-1\Big].$
	\\
	\\
	(b)  If $g$ is even, then  $${\rm avec}(G)\leq   \frac{3g}{4}\Big( \frac{n-L_2}{2L_1}\Big)\Big(1+ \frac{L_2-L_1}{3n} \Big) +\frac{21g}{8}-2,$$ where 
	$ L_1= \frac{(\delta-1)^{\frac{g}{2}}-1}{\delta-2}$ and $L_2 = \Delta + \frac{\Delta-1}{\delta -2} \Big[ (\delta-1)^{\frac{g-2}{2}} -(\delta-1)\Big] $.
\end{thm}
\noindent
\textbf{Proof.}
(a)\,\, Let $v_1$ be a vertex of degree $\Delta$ and let $A = \{v_1, a_1,a_2, \ldots, a_r\}$ be a maximal $(g-1)$-packing  constructed such that each vertices of $G$ not in $A$ is withing distance $g-1$ of $A$.  Let $T$ be a spanning tree obtained in the same manner as in the proof of Theorem \ref{theoremgirth} (a). Then, $deg_T(v_1) = deg_G(v_1)$. For every vertex $x$ of $T$ let $x_A$ be a vertex in $A$ closest to $u$ in $T$ and let $c:V(T)\rightarrow \mathbb{R}^+$ be the weight function defined by 
$$c(u) = \Big|\{ x \in V(T)|  x_A=u\}\Big|, \, \text{for each}\,\, u\in A.$$ 
Then,
\begin{enumerate}
	\item[(i.)] $c(u) = 0$ if $u \not \in A$.
	\item[(ii.)] For each $u \in A -\{v_1\}$, \begin{equation}
	\begin{aligned}
	c(u) &\geq |{N}_G^{(g-1)/2}[u]|\\
	&\geq 1+ \delta +\delta (\delta-1) + \delta (\delta-1)^2+ \ldots + \delta(\delta-1)^{(g-3)/2}.\\
	&= 1+ \frac{\delta}{\delta-2}\big[ (\delta-1)^{(g-1)/2}-1\big] = K_1.
	\end{aligned}
	\end{equation}
	\item[(iii.)] For $v_1 \in A$, 
	\begin{equation}
	\begin{aligned}
	c(v_1) &\geq |{N}_G^{(g-1)/2}[v_1]|\\
	&\geq 1+ \Delta +\Delta (\delta-1) + \Delta (\delta-1)^2+ \ldots + \Delta(\delta-1)^{(g-3)/2}.\\
	&= 1+ \frac{\Delta}{\delta-2}\big[ (\delta-1)^{(g-1)/2}-1\big] = K_2.\\
	\end{aligned}
	\end{equation}
\end{enumerate}
Since $n= \sum_{u\in A} c(u)$, it follows that $n\geq K_2 + (|A| - 1)K_1 $ and rearranging yields,
\begin{equation} 
|A| \leq \frac{n-K_2}{K_1}+1.
\end{equation}
By the construction of $A$, it is easy to see that inequations \eqref{eqn2p} and \eqref{eqn2q} are true for Theorem \ref{theoremgirth2} (a), as well as the fact that $T^g[A]$ is connected.
Let $c'(u)$ be a new weight function which satisfies $c'(u) \geq 1$ for all $u \in A$. Define $c'(u)=\frac{c(u)}{K_1}$ for $u\in A-\{v_1\}$ and $c'(v_1)=\frac{c(v_1)-K_2+K_1}{K_1}$. Thus, we have the total weight $N$ of $c'$ as
\[ N = \sum_{u \in A}c'(u)= \frac{n-K_2}{K_1}+1,\] implying that $|A| \leq N$.
We simplify ${\rm avec}_{c'}(T^g[A])$ to yield
\begin{eqnarray*}
	{\rm avec}_{c'}(T^g[A]) 
	&=& \frac{EX_{c'}(T^g[A])}{N} \\
	&=&  \frac{EX_{c}(T^g[A])-e(v_1)(K_2-K_1)}{n-K_2-K_1}, 
\end{eqnarray*}
which in conjunction with $EX(T^g[A]) = n \,\,{\rm avec}_c(T^g[A])$ and by rearranging becomes 
\begin{equation}\label{eqn9}
{\rm avec}_{c}(T^g[A]) = \frac{n-K_2+K_1}{n} {\rm avec}_{c'}(T^g[A]) - \frac{K_2-K_1}{n} e_{T^g[A]}(v_1).
\end{equation}
Since the order of $T^g[A]$ is $|A|$ and since $|A| \leq \frac{n-K_2}{K_1} +1$, we have
\begin{eqnarray*}
	e_{T^g[A]}(v_1) 
	&\leq &  {\rm diam} (T^g[A]) \leq |A|-1\\
	&\leq & \frac{n-K_2}{K_1}.
\end{eqnarray*}
By Proposition \ref{prop2} and in conjunction with $\lceil N \rceil < \frac{n-K_2}{K_1}+2$, we have
\begin{eqnarray*}
	{\rm avec}_{c'}(T^g[A]) 
	&\leq & \frac{3\lceil N \rceil}{4} -\frac{1}{2} < \frac{3}{4}\Big( \frac{n-K_2}{K_1} +2 \Big) -\frac{1}{2}= \frac{3}{4}\Big(  \frac{n-K_2}{K_1}\Big) +1. \\
\end{eqnarray*}
By substituting the bounds for $e_{T^g[A]}(v_1)$ and ${\rm avec}_{c'}(T^g[A])$ in \eqref{eqn9}, we have
\begin{eqnarray*} 
	{\rm avec}_c(T^g[A]) 
	&\leq &  \Big(\frac{n-K_2+K_1}{n}\Big)\Big(\frac{3}{4} \frac{n-K_2}{K_1} + 1 \Big) + \frac{K_2-K_1}{n}\Big(\frac{n-K_2}{K_1}\Big)  \\
	&\leq &  \frac{3}{4} \Big( \frac{n-K_2}{K_1}\Big) +\frac{3}{4}\Big( \frac{n-K_2}{K_1}\Big)\Big( \frac{K_2-K_1}{3n}\Big) + 1\\
	&=& \frac{3}{4} \Big( \frac{n-K_2}{K_1}\Big) \Big(1 + \frac{K_2-K_1}{3n}\Big) +1.
\end{eqnarray*}
Since inequalities \eqref{eqn2p}, \eqref{eqn2q} and \eqref{eqn2r} are also true for Theorem \ref{theoremgirth2}, we have
\begin{eqnarray*}
	{\rm avec}(T) 
	&\leq & {\rm avec}_c(T) +g-1\\
	&\leq & g\,\, {\rm avec}_c(T^g[A]) +2(g-1)\\
	&\leq & g\Big[ \frac{3}{4} \Big( \frac{n-K_2}{K_1}\Big) \Big(1 + \frac{K_2-K_1}{3n}\Big) +1 \Big] + 2(g-1)\\
	&=& \frac{3g}{4} \Big( \frac{n-K_2}{K_1}\Big) \Big(1 + \frac{K_2-K_1}{3n}\Big)  +3g -2.
\end{eqnarray*}
The results follows since ${\rm avec}(G) \leq {\rm avec}(T)$. 
\\
\\
(b)\,\, For $g$ even, let $v_1$ be a vertex of degree $\Delta$ and let $e_1$ be an edge incident with $v_1$. Let $M=\{e_1,e_2,\ldots, e_r\}$ be a maximum matching such that each of the edges not in $M$ is within distance $g-2$ of $M$. Let $T$ be a spanning tree of $G$ obtained as in the proof of Theorem \ref{theoremgirth} $(b)$. Define $V(M)$ as the set of vertices incident with an edge in $M$. For every vertex $u\in V(T)$, let $u_M$ be a vertex in $V(M)$ closest to $u$ in $T$. We define a weight function $c:V(T)\rightarrow \mathbb{R}^+$ by 
$$c(u) = \Big|\{ x \in V(T)|  x_M=u\}\Big|, \, \text{for}\,\, u\in V(T).$$  Since $deg_T(v_1)=deg_G(v_1)$, $T$ has the same maximum degree in $G$. Furthermore, since $g$ is even, the least possible value of $g$ is 4. Thus, $G$ is triangle-free and so, no two incident verties of any edge in $M$ have a common neighbour. It follows that 
\begin{enumerate}
	\item[(i.)] $c(u) = 0$ if $u \not \in M$.
	\item[(ii.)] For each $u \in M -\{v_1\}$, \begin{equation}
	\begin{aligned}
	c(u) &\geq |{N}_G^{(g-2)/2}[u]|\\
	&\geq 1+ (\delta-1) + (\delta-1)^2+ \ldots + (\delta-1)^{(g-2)/2}\\
	&= \frac{(\delta-1)^{g/2}-1}{\delta -2} =L_1.
	\end{aligned}
	\end{equation}
	\item[(iii.)] For $v_1 \in M$, 
	\begin{equation}
	\begin{aligned}
	c(v_1) &\geq |{N}_G^{(g-1)/2}[v_1]|\\
	&\geq 1+ (\Delta-1) + (\Delta-1) (\delta-1) + (\Delta-1) (\delta-1)^2+ \ldots + (\Delta-1)(\delta-1)^{(g-4)/2} \\
	&= \Delta + \frac{\Delta-1}{\delta -2} \Big[ (\delta-1)^{\frac{g-2}{2}} -(\delta-1)\Big] = L_2.
	\end{aligned}
	\end{equation}
\end{enumerate} 
Let $\bar{c}$ be the weight function on the line graph $L(T)$ with $V(L)=E(T)$ defined by:
\[
\bar{c}(uv) =
\begin{cases}
c(u) +c(v) & \text{if } uv \in M, \\
0
& \text{if } uv \not \in M.
\end{cases}
\] 
Since $e_1$ is incident with $v_1$ in $T$, $\bar{c}(e_1) \geq L_1+L_2$ and $\bar{c}(e) \geq 2L_1$ for all $e\in M -\{e_1\}$. It follows that since $n= \sum_{e \in M}\bar{c}(e)$, 
\[n\geq L_2 +L_1 + 2L_1(|M|-1),\] and rearranging yields
\begin{equation}
|M| \leq \frac{n-L_2+L_1}{2L_1}.
\end{equation}
As illustrated in Theorem \ref{theoremgirth} (b), inequalities \eqref{2.20}, \eqref{2.21} and \eqref{2.22} also hold for Theorem \ref{theoremgirth2} (b). To apply Proposition \ref{prop2} to ${\rm avec_{\bar{c}}}(T^g[M])$, we normalise the weight function $\bar{c}$.
Let $\bar{c'}$ be a new weight function defined by $\bar{c'}(e)\geq 1$ satisfying $\bar{c'}(e)= \frac{\bar{c}(e)}{L_1}\geq 1$ and $\bar{c'}(e_1) = \frac{\bar{c}(e_1)-L_2+L_1}{2L_1}\geq 1$ for all $e \in E(G)$. Thus,
\[N'=\sum_{e \in M} \bar{c'}(e) = \frac{n-L_2}{2L_1} +\frac{1}{2}.\] 
Then, as shown in Theorem \ref{theoremgirth2} (a) (see \ref{eqn9}), 
\begin{equation}\label{eqn*}
{\rm avec}_{{\bar{c}}}(L^g[M]) = \frac{n-L_2+L_1}{n} {\rm avec}_{\bar{c'}}(L^g[M])- \frac{L_2-L_1}{n} e_{L^g[M]}(v_1).
\end{equation}
Since the order of $L^g[M]$ is $|M|$ and since $|M| \leq \frac{n-L_2}{2L_1} +\frac{1}{2}$, we have
\begin{eqnarray*}
	e_{T^g[A]}(v_1) 
	&\leq &  {\rm diam} (L^g[M]) \leq |M|-1\\
	&\leq & \frac{n-L_2}{2L_1} -\frac{1}{2}.
\end{eqnarray*} 
Letting $\lceil N' \rceil < N' + 1 = \frac{n-L_2}{2L_1} +\frac{3}{2}$, we have by Proposition \ref{prop2}, 
\[{\rm avec}_{\bar{c'}}(L^g[M]) \leq \frac{3\lceil N' \rceil}{4}  -\frac{1}{2} <\frac{3}{4} \Big( \frac{n-L_2}{2L_1}+ \frac{3}{2}   \Big) - \frac{1}{2}= \frac{3}{4} \Big( \frac{n-L_2}{2L_1}\Big) + \frac{5}{8}. \]
By substituting inequalities for ${\rm avec}_{\bar{c'}}(L^g[M])$ and $e_{L^g[M]}(v_1)$ in \eqref{eqn*} and rearranging yields
\begin{eqnarray*}
	{\rm avec}_{\bar{c}}(L^g[M]) &\leq & \Big(\frac{n-L_2+L_1}{n}\Big)\Big(\frac{3}{4}\frac{n-L_2}{2L_1}+\frac{5}{8}\Big) +\frac{L_2-L_1}{n}\Big(\frac{n-L_2}{2L_1}-\frac{1}{2}\Big).\\
	&\leq & \frac{3}{4} \Big( \frac{n-L_2}{2L_1}\Big)\Big[ 1+\Big( \frac{L_2-L_1}{3n}\Big)\Big] +\frac{5}{8}.
\end{eqnarray*}
Since the inequalities \eqref{2.20}, \eqref{2.21} and \eqref{2.22} also holds for Theorem \ref{theoremgirth2} (b), we have
\begin{equation*}
\begin{aligned}
{\rm avec}(T) &\leq  {\rm avec}_c(T) + g-2\\
& \leq g\, {\rm avec}_{\bar{c}}(L) +g-1\\
& \leq g\, {\rm avec}_{\bar{c}}(L^g[M]) +2(g-1)\\
& \leq g\Big[  \frac{3}{4} \Big( \frac{n-L_2}{2L_1}\Big)\Big[ 1+\Big( \frac{L_2-L_1}{3n}\Big)\Big] +\frac{5}{8} \Big] +2g-2\\
&= \frac{3g}{4}\Big( \frac{n-L_2}{2L_1}\Big)\Big[ 1+\Big( \frac{L_2-L_1}{3n}\Big)\Big] +\frac{21g}{8}-2.
\end{aligned}
\end{equation*}
Since deleting edges does not decrease the average eccentricity, ${\rm avec}(G)\leq {\rm avec}(T)$, and thus Theorem \ref{theoremgirth2} follows as desired.
\hfill $\Box$
\\
\\
It is easy to verify that the results are equivalent for cases $g=3, 4$ and 5, as proved in \cite{DanOsa2019}.
\begin{enumerate}
	\item[(a)] If $g=3$, $G$ represents any connected graph of given $n$, $\delta \geq 3$ and $\Delta$. Then $K_1 = \delta +1$ and $K_2= \Delta +1$, and so by a simple calculation, we have ${\rm avec}(G) \leq \frac{9}{4} \frac{n-\Delta-1}{\delta+1} \big(  1+ \frac{\Delta -\delta}{3n}\big) +7$. Note that this is equivalent to the bound \eqref{eqn6}. 
	\item [(b)] If $g=4$, $G$ is triangle-free. Then $L_1 = \delta$ and $L_2 = \Delta$, and so ${\rm avec} (G) \leq 3 \big( \frac{n-\Delta}{2\delta} \big) \big(1+ \frac{\Delta -\delta}{3n} \big) + \frac{21}{2}$, which differs from \eqref{eqn7} only by the additive constant 1.
	\item [(c)]If $g=5$, then $G$ is a $(C_3,C_4)$-free graph with $K_1= \delta^2+1$ and $K_2= \Delta \delta +1$.  Therefore, ${\rm avec}(G) \leq \frac{15}{4} \big(\frac{n-\delta \Delta -1}{\delta^2 +1} \big) \big(1+ \frac{\delta(\Delta -\delta)}{3n}  \big) +13 $. This differs from \eqref{eqn8} for $C_4$-free graphs by some additive constant. 
\end{enumerate}
Theorem \ref{theoremgirth2} generalises Theorem \ref{theoremgirth} in the sense that setting $\Delta$ equal to $\delta$ yields $K_1=K_2=K$ and $L_1=L_2=L$.

\subsection*{Acknowledgement}
The author, most sincerely, thank Dr. Peter Dankelmann for his immense support and mentorship during the course of my doctoral studies at the University of Johannesburg and for his helpful suggestions on this paper.  

%


\begin{thebibliography}{99}


 
%
%


%
%
%
%
%



\bibitem{PAli} P.\ Ali,  
The Steiner diameter of a graph with prescribed girth.
Discrete \ Math.\ {\bf 313} (2016), 1322-1326. 

\bibitem{AliDanMorMukSwaVet2018} P.\ Ali, P.\ Dankelmann, M.J.\ Morgan, 
S.\ Mukwembi, T.\ Vetr\'\i k, 
The average eccentricity, spanning trees of plane graphs, size
and order. Utilitas Math.\ {\bf 107} (2018), 37-49. 

\bibitem{AlexDan2020} A. Alochukwu, P. Dankelmann, 
Upper bounds on the average eccentricity of graphs of girth 6 and $(C_4,C_5)$-free graphs, (submitted) https://arxiv.org/pdf/2004.14490.pdf.   	

\bibitem{Aouchiche} M.\ Aouchiche, G.\ Caporossi, P.\ Hansen, 
M.\ Laffay, 
AutoGraphiX: a survey, Electron. Notes Discrete Math.\ {\bf 22} (2005), 515-520.  

\bibitem{BucHar} F.\ Buckley, F.\ Harary,
Distance in Graphs.
Addisson-Wesley, Redwood City, California (1990).    

\bibitem{DanEnt2000} P.\ Dankelmann, R.\ Entringer,
Average distance, minimum degree and spanning trees.
J.\ Graph Theory {\bf 33} no.\ 1 (2000), 1-13.  

\bibitem{DanGodSwa2004} P.\ Dankelmann,  W.\ Goddard, C.S.\ Swart, 
The average eccentricity of a graph and its subgraphs.  
Util.\ Math.\ {\bf 41} (2004), 41-51.

\bibitem{DanMuk2014} P.\ Dankelmann and S.\ Mukwembi, 
Upper bounds on the average eccentricity. 
Discrete Appl.\ Math.\ {\bf 167} (2014), 72-79. 

\bibitem{DanOsa-manu} P.\ Dankelmann, F.J.\ Osaye, 
Average eccentricity, $k$-packing and $k$-domination in graphs. Discrete \ Math. \ {\textbf{342}} (2019), 1261- 1274.

\bibitem{DanOsa2019} P.\ Dankelmann, F.J.\ Osaye, 
Average eccentricity, minimum degree and maximum degree in graphs. J. Combin. Optim. (2020). 


\bibitem{DanMukOsaRod} P.\ Dankelmann, F.J.\ Osaye, S. \ Mukwembi, B. \ Rodrigues, 
Upper Bounds on the	average eccentricity of $K_3$-free and $C_4$-free graphs.
Discrete Appl. Math., \textbf{270} (2019), 106-114.

\bibitem{DuIli2013} Z.\ Du, A.\ Ili\u{c}, 
On AGX conjectures regarding average eccentricity.        
MATCH Commun.\ Math.\ Comput.\ Chem.\ {\bf 69} (2013), 597-609.  


\bibitem{ErdPacPolTuz1989} P.\ Erd\"{o}s, J.\ Pach, R.\ Pollack, and Z.\  Tuza, 
Radius, diameter, and minimum degree, 
J.\ Combin.\ Theory B {\bf 47} (1989), 73-79.                     

\bibitem{Ili2012} A.\ Ili\u{c},            
On the extremal properties of the average eccentricity.
Computers and Mathematics with Applications {\bf 64} no.\ 9 (2012), 2877-2885.  

\bibitem{SmiSzeWan2016} H.\ Smith, L.A.\ Sz\'{e}kely, H.\ Wang,
Eccentricity sum in trees.
Discrete Appl.\ Math.\ {\bf 207} (2016), 120-131.  

\bibitem{TanZho2012} Y.\ Tang, B.\ Zhou,
On average eccentricity. 
MATCH Commun.\ Math.\ Comput.\ Chem.\ {\bf 67} (2012), 405-423.  
\end{thebibliography}
\end{document}